\documentclass[11pt]{article}
\usepackage{enumerate}

\usepackage{amssymb,a4wide,latexsym,makeidx,epsfig,fleqn}
\usepackage{amsthm}
\usepackage{amsmath}
\usepackage{enumerate}
\usepackage{graphicx}
\usepackage{subfigure}
\usepackage{float}
\usepackage{caption}
\usepackage[colorlinks, linkcolor=black, anchorcolor=black, citecolor=blue]{hyperref}%
\newtheorem{theorem}{Theorem}[section]

\newtheorem{definition}[theorem]{Definition}
\newtheorem{lemma}[theorem]{Lemma}

\newtheorem{problem}[theorem]{Problem}

\begin{document}
\textwidth 150mm \textheight 225mm
\title{Maximizing the signless Laplacian spectral radius of some theta graphs
\thanks{ Supported by the National Natural Science Foundation of China (No. 12271439).}}
\author{{Yuxiang Liu$^{a,b}$, Ligong Wang$^{a,b,}$\footnote{Corresponding author.}}\\
{\small $^{a}$  School of Mathematics and Statistics, Northwestern Polytechnical University,}\\{\small  Xi'an, Shaanxi 710129, P.R. China.}\\ {\small $^{b}$ Xi'an-Budapest Joint Research Center for Combinatorics, Northwestern Polytechnical University,}\\{\small  Xi'an, Shaanxi 710129, P.R. China.}\\
{\small E-mail: yxliumath@163.com, lgwangmath@163.com}}
\date{}
\maketitle
\begin{center}
\begin{minipage}{135mm}
\vskip 0.3cm
\begin{center}
{\small {\bf Abstract}}
\end{center}
Let $Q(G)=D(G)+A(G)$ be the signless Laplacian matrix of a simple graph $G$, where $D(G)$ and $A(G)$ are the degree diagonal matrix and the adjacency matrix of $G$, respectively. The largest eigenvalue of $Q(G)$, denoted by $q(G)$, is called the signless Laplacian spectral radius of $G$. Let $\theta(l_{1},l_{2},l_{3})$ denote the theta graph which consists of two vertices connected by three internally disjoint paths with length $l_{1}$, $l_{2}$ and $l_{3}$. Let $F_{n}$ be the friendship graph consisting of $\frac{n-1}{2}$ triangles which intersect in exactly one common vertex for odd $n\geq3$ and obtained by hanging an edge to the center of $F_{n-1}$ for even $n\geq4$. Let $S_{n,k}$ denote the graph obtained by joining each vertex of $K_{k}$ to $n-k$ isolated vertices. Let $S_{n,k}^{+}$ denote the graph obtained by adding an edge to the two isolated vertices of $S_{n,k}$. In this paper, firstly, we show that if $G$ is $\theta(1,2,2)$-free, then $q(G)\leq q(F_{n})$, unless $G\cong F_{n}$. Secondly, we show that if $G$ is $\theta(1,2,3)$-free, then $q(G)\leq q(S_{n,2})$, unless $G\cong S_{n,2}$. Finally, we show that if $G$ is $\{\theta(1,2,2),F_{5}\}$-free, then $q(G)\leq q(S_{n,1}^{+})$, unless $G\cong S_{n,1}^{+}$.

\end{minipage}
\end{center}
\section{Introduction}
In this paper, we consider only simple and finite graphs and we follow the traditional notation and terminology \cite{BoMu1}. For a graph $G=(V(G), E(G))$, we use $n=|V(G)|$ and $m=e(G)=|E(G)|$ to denote the order and the size of $G$, respectively. For a vertex $u\in V(G)$, let $N_{G}(u)$ be the neighborhood of $u$ and $N_{G}[u]=N_{G}(u)\cup \{u\}$. Let $d_{G}(u)$ be the degree of $u$ and $\Delta(G)=\Delta$ be the maximum degree of $G$. For the sake of simplicity, we omit all the subscripts defined here if $G$ is clear from the context, for example, $N(u)$ and $d(u)$.

Denote by $A(G)$ the adjacency matrix of $G$. The largest eigenvalue of $A(G)$, denoted by $\lambda(G)$, is named as the adjacency spectral radius of $G$. Denote by $D(G)$ the diagonal degree matrix. The matrix $Q(G)=D(G)+A(G)$ is the signless Laplacian matrix of $G$. The largest eigenvalue of $Q(G)$, denoted by $q(G)$, is named as the signless Laplacian spectral radius of $G$. From Perron-Forbenius theorem, for a connected graph $G$ there exists a positive unit eigenvector corresponding to $q(G)$, which is called the Perron vector of $Q(G)$.

Let $\theta(l_{1},l_{2},l_{3})$ denote the theta graph which consists of two vertices connected by three internally disjoint paths with length $l_{1}$, $l_{2}$ and $l_{3}$. A generalized theta graph $\theta(l_{1},l_{2},\cdots, l_{t})$ is the graph obtained by connecting two vertices with $t$ internally disjoint paths of lengths $l_{1}, l_{2}, \cdots, l_{t}$, where $l_{1}\leq l_{2}\leq\cdots \leq l_{t}$ and $l_{2}\geq 2$. Let $F_{n}$ be the friendship graph consisting of $\frac{n-1}{2}$ triangles which intersect in exactly one common vertex for odd $n\geq3$ and obtained by hanging an edge to the center of $F_{n-1}$ for even $n\geq4$. Let $S_{n,k}$ denote the graph obtained by joining each vertex of $K_{k}$ to $n-k$ isolated vertices. Let $S_{n,k}^{+}$ denote the graph obtained by adding a edge to the two isolated vertices of $S_{n,k}$.

A graph is said to be $\mathcal{F}$-free, if it does not contain any graph in $\mathcal{F}$ as a subgraph. In particular, if $\mathcal{F}=\{F\}$, we also say that $G$ is $F$-free. Let $\mathcal{G}(n,\mathcal{F})$ denote the family of $\mathcal{F}$-free graphs with $n$ vertices and without isolated vertices. In particular, if $\mathcal{F}=\{F\}$, then we write $\mathcal{G}(n, \mathcal{F})$ for $\mathcal{G}(n,F)$. The $Tur\acute{a}n$ number, denoted by $ex(n,\mathcal{F})$, is the maximal number of edges in a $\mathcal{F}$-free graph of order $n$. To determine the exact value of $ex(n,\mathcal{F})$ is a central problem of extremal graph theory. For more details on this topic, the readers may referred to \cite{Bol,Nik7,Tur} and the references therein. This problem was considered about theta graphs, see e.g., $F=\theta(l_{1},l_{2},\cdots, l_{t})$, where $l_{1}=l_{2}=\cdots= l_{t}=k$ \cite{BuTa,FaSi}; $F=\theta(1,2,4)$ \cite{LiSW1}; $F=\theta(l_{1},l_{2},\cdots, l_{t})$, where the lengths of all $l_{i}, i\in \{1,2,\cdots,t\}$ have the same length and at most one $l_{i}, i\in \{1,2,\cdots,t\}$ has length 1 and $F=\theta(3,5,5)$ \cite{LiYa}; $F=\theta(4,4,4)$ \cite{VeWi}; $F=\theta(p,q,r)$ \cite{ZhFS}, where $p\leq q\leq r$ and $q\geq2$.

In 2010, Nikiforov \cite{Nik1} proposed Spectral Tur\'{a}n type problem about determining the maximum spectral radius of an $F$-free graph of order $n$. Over the past decade, much attention has been paid to the topic. For more details, one may refer to \cite{ChZ,LiLF} and the references therein. This problem was considered about theta graphs, see e.g., $F=\theta(1,2,r)$, where $n\geq10r$ if $r$ is odd and $n\geq 7r$ if $r$ is even \cite{ZhLi}; $F=\theta(1,2,3)$ and $F=\theta(1,2,4)$ \cite{LiSW1}. In 1985, Brualdi and Hoffman \cite{BrHo} proposed the problem about determining the maximum adjacency spectral radius of an $F$-free graph of size $m$. This problem was considered about theta graphs, see e.g., $\mathcal{F}=\{\theta(1,2,2),\theta(1,2,3)\}$ \cite{SLW, ZhLS}; $F=\theta(1,2,3)$ \cite{FaY, SLW}; $F=\theta(1,2,4)$ \cite{SLW}. In 2013, de Freitas, Nikiforov and Patuzzi \cite{FNP1} proposed the spectral extremal problem of signless Laplacian matrix of graph as follows.

\noindent\begin{problem}
 What is the maximum signless Laplacian spectral radius of an $F$-free graph of order $n$?
\end{problem}
Problem 1.1 was considered for various forbidding subgraphs such as $kP_{2}$ \cite{Yu}; paths \cite{NiYu1}; $C_{4}$ and $C_{5}$ \cite{FNP1}; odd cycles \cite{Yuan}; even cycles \cite{NiYu2}; $K_{s,t}$ \cite{FNP2}, linear forests \cite{CLZ}; friendship graphs \cite{Zhao}; flowers \cite{ChZh}; fan \cite{WaZh}. For the related results, one may refer to \cite{CvS,CvS1,CvS2} and the references therein. To our knowledge, we can not find out the results about the forbidding theta graph and we get the results as follows.

\noindent\begin{theorem}\label{th:ch-1.2}
Let $G$ be a $\theta(1,2,2)$-free graph with order $n\geq6$, then $q(G)\leq q(F_{n})$, unless $G\cong F_{n}$.
\end{theorem}

\noindent\begin{theorem}\label{th:ch-1.3}
Let $G$ be a $\theta(1,2,3)$-free graph with order $n\geq6$, then $q(G)\leq q(S_{n,2})$, unless $G\cong S_{n,2}$.
\end{theorem}

\noindent\begin{theorem}\label{th:ch-1.4}
If $G$ is a $\{\theta(1,2,2),F_{5}\}$-free graph with order $n\geq6$, then $q(G)\leq q(S_{n,1}^{+})$, with equality holds if and only if $G\cong S_{n,1}^{+}$.
\end{theorem}
\section{Preliminary}

For convenience, We need introduce some notations. For a graph $G$ and a subset $S\subseteq V(G)$, let $G[S]$ denote the subgraph of $G$ induced by $S$. For a vertex $u$ and a vertex subset $S$, $N_{S}(u)=N(u)\cap S$ and $d_{S}(u)=|N_{S}(u)|$. Let $N^{2}(u)$ be the set of vertices of distance two to $u$. Let $N_{i}(u)=\{v| v\in N(u), d_{N(u)}(v)=i\}$. Let $X$ and $Y$ be two disjoint sets of vertices of $G$. Denote by $e(X,Y)$ the number of edges joining every vertex in $X$ to every vertex in $Y$. Given two disjoint graphs $G$ and $H$, let $G\vee H$ denote the graph obtained by joining every vertex of $G$ to every vertex of $H$.

Throughout the paper, let $G^{\ast}$ denote an extremal graph with maximal signless Laplacian spectral radius in $\mathcal{G}(n,\mathcal{F})$. Let $N[u]=N(u)\cup \{u\}$, $W=V(G^{\ast})\backslash N[u]$ and $N_{i}(u)=\{v| v\in N(u), d_{N(u)}(v)=i, i\in\{0,1\}\}$. Let $u\in V(G^{\ast})$ be the vertex such that
\begin{equation}\label{eq:ch-1}
d(u)+\frac{1}{d(u)}\sum_{w\in N(u)}d(w)=max\{d(v)+\frac{1}{d(v)}\sum_{w\in N(v)}d(w): v\in V(G^{\ast})\}.
\end{equation}

\noindent\begin{definition}\label{le:ch-2.1}{\rm(}$\cite{CvRS}${\rm)} Given a graph $G$, the vertex partition $\Pi$: $V(G)=V_{1}\cup V_{2} \cup \ldots \cup V_{k}$ is said to be an equitable partition if, for each $u\in V_{i}$, $|V_{j}\cap N(u)|=b_{ij}$ is a constant depending only on $i,j$ ($1\leq i,j\leq k$). The matrix $B_{\Pi}=(b_{ij})$ is called the quotient matrix of $G$ with respect to $\Pi$.
\end{definition}

\noindent\begin{lemma}\label{le:ch-2.2}{\rm(}$\cite{CvRS}${\rm)} Let $\Pi$: $V(G)=V_{1}\cup V_{2} \ldots \cup V_{k}$ be an equitable partition of $G$ with quotient matrix $B_{\Pi}$. Then $det(xI-B_{\Pi}) \mid det(xI-A(G))$. Furthermore, the largest eigenvalue of $B_{\Pi}$ is just the spectral radius of $G$.
\end{lemma}

\noindent\begin{lemma}\label{le:ch-2.3}{\rm(}$\cite{ErGa}${\rm)}
If $G$ is a $P_{k+1}$-free graph of order $n$, then $e(G)\leq\frac{(k-1)n}{2}$, with equality if and only if $G$ is a union of disjoint copies of $K_{k}$.
\end{lemma}

\noindent\begin{lemma}\label{le:ch-2.4}{\rm(}$\cite{FeYu}${\rm)}
For any graph $G$, we have
$$
q(G)\leq max_{u\in V(G)}\{d(u)+\frac{1}{d(u)}\sum_{w\in N(u)}d(w)\},
$$
with equality if and only if $G$ is either a semi-regular bipartite graph or a regular graph.
\end{lemma}

\noindent\begin{lemma}\label{le:ch-2.5}{\rm(}$\cite{Das}${\rm)}
Let $G$ be a simple graph with $n$ vertices and $m$ edges. We have

$$
max\{d(u)+\frac{1}{d(u)}\sum_{w\in N(u)}d(w)\}\leq \frac{2m}{n-1}+n-2,
$$
with equality if and only if $G$ is either complete graph, or a star, or a complete
graph with one isolated vertex.
\end{lemma}

\noindent\begin{lemma}\label{le:ch-2.6}{\rm(}$\cite{ChTa}${\rm)}
If $3\leq k\leq n-3$, then $H_{n,k}$ the graph obtained from the star $K_{1,n-1}$ by joining a vertex of degree 1 to $k+1$ other vertices of degree 1, is the unique connected graph that maximizes the signless Laplacian spectral radius over all connected graphs with $n$ vertices and $n+k$ edges.
\end{lemma}

\noindent\begin{lemma}\label{le:ch-2.7}{\rm(}$\cite{CvS}${\rm)}
If $H$ is the subgraph of a connected graph $G$, then $q(H)\leq q(G)$. Particularly, if $H$ is proper, then $q(H)<q(G)$.
\end{lemma}

\noindent\begin{lemma}\label{le:ch-2.8}{\rm(}$\cite{FNP1}${\rm)} {\rm(}$i${\rm)} If $n\geq3$ is odd, then $q(F_{n})=\frac{n+2+\sqrt{(n-2)^{2}+8}}{2}$ and satisfies
$$
 n+\frac{2}{n-1}<q(F_{n})<n+\frac{2}{n-2};
$$
{\rm(}$ii${\rm)} If $n\geq4$ is even, then $q(F_{n})$ is the largest root of equation $x^{3}-(n+3)x^{2}+3nx-2n+4=0$, and satisfies
$$
n+\frac{2}{n}<q(F_{n})<n+\frac{2}{n-1}.
$$
\end{lemma}

\noindent\begin{lemma}\label{le:ch-2.9}{\rm(}$\cite{FNP1}${\rm)}
If $n\geq4$, then
$$
q(S_{n,2})=\frac{n+2+\sqrt{n^{2}+4n-12}}{2}>n+2-\frac{4}{n+1}.
$$
\end{lemma}

\section{Proof of Theorem 1.2}

\noindent {\bf Proof.} Since $G^{\ast}$ is $\theta(1,2,2)$-free, we obtain $G^{\ast}[N(u)]$ is $P_{3}$-free, and hence $N(u)=N_{0}(u)\cup N_{1}(u)$. Since $F_{n}$ is $\theta(1,2,2)$-free, we have $$n+\frac{2}{n-1}<q(F_{n})\leq q(G^{\ast})$$ for odd $n$ and $$n+\frac{2}{n}<q(F_{n})\leq q(G^{\ast})$$ for even $n$ from Lemma \ref{le:ch-2.8}. Let $u\in V(G^{\ast})$ be the vertex such that equality \eqref{eq:ch-1} holding. Obviously, for any $w\in W$, $|N_{N_{1}(u)}(w)|\leq\frac{|N_{1}(u)|}{2}$. Suppose that $G^{\ast}[N(u)]$ contains $tK_{2}$, where $0\leq t\leq\lfloor\frac{d(u)}{2}\rfloor$, we have
$$
\begin{aligned}
e(N(u),N^{2}(u))&\leq (n-1-d(u))d(u)-(n-1-d(u))t\\&=(n-1-d(u))(d(u)-t).
\end{aligned}
$$
Furthermore,
\begin{equation}\label{eq:ch-2}
\begin{split}
q(G^{\ast})&\leq d(u)+\frac{1}{d(u)}\sum_{w\in N(u)}d(w)\\&= d(u)+\frac{d(u)+2t+(d(u)-t)(n-1-d(u))}{d(u)}\\&= n-t(\frac{n-3-d(u)}{d(u)}).
\end{split}
\end{equation}

{\bf Case 1.} If $\Delta=n-1$, then there exists a vertex $v$ such that $d(v)=n-1$. Recall that $N(v)=N_{0}(v)\cup N_{1}(v)$. Furthermore, we obtain $G^{\ast}$ is a subgraph of $F_{n}$ and $q(G^{\ast})<q(F_{n})$ from Lemma \ref{le:ch-2.7}, unless $G^{\ast}\cong F_{n}$, as desired.

{\bf Case 2.} If $\Delta\leq n-2$, then we consider two subcases as follows.

{\bf Subcase 2.1.} If $d(u)\leq n-3$, by Lemma \ref{le:ch-2.8} and \eqref{eq:ch-2}, then both
$$
\begin{aligned}
 n+\frac{2}{n}<q(G^{\ast})\leq n-t(\frac{n-3-d(u)}{d(u)})\leq n,
\end{aligned}
$$
and
$$
\begin{aligned}
n+\frac{2}{n-1}<q(G^{\ast})\leq n-t(\frac{n-3-d(u)}{d(u)})\leq n,
\end{aligned}
$$
a contradiction.

\begin{figure}[H]
 \begin{center}
   \begin{minipage}{0.5\textwidth}
    \includegraphics[height=3cm,width=5cm]{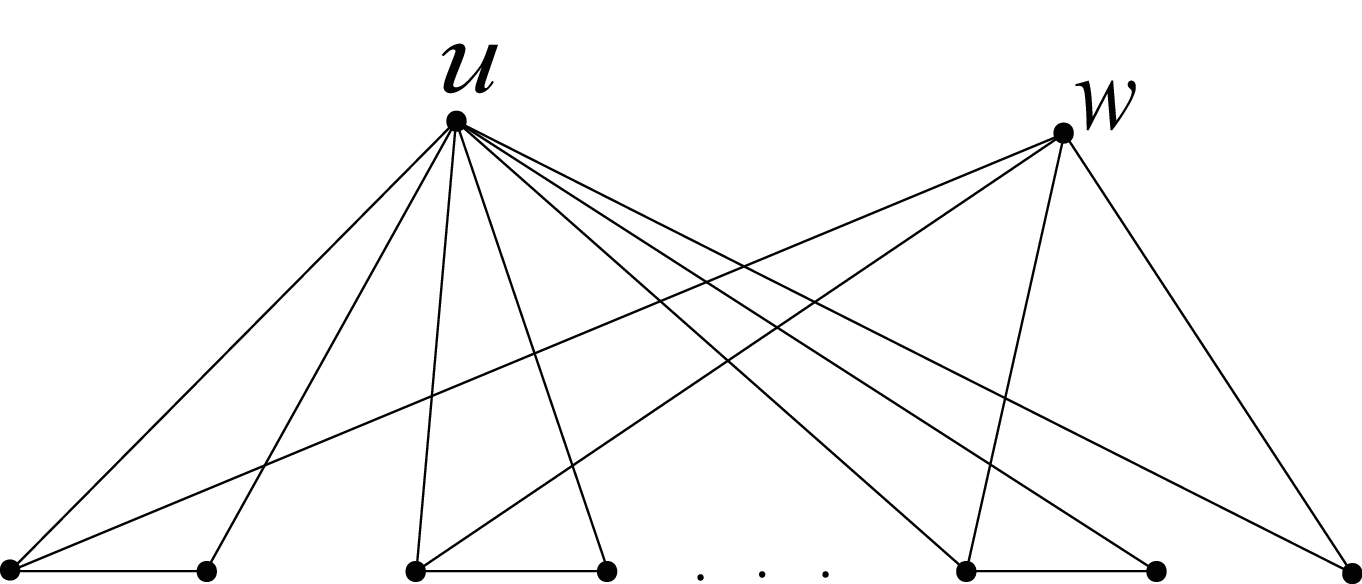}
     \caption{The graph $G_{1}$.}
      \label{1.}
   \end{minipage}
  \begin{minipage}{0.5\textwidth}
   \includegraphics[height=3cm,width=5cm]{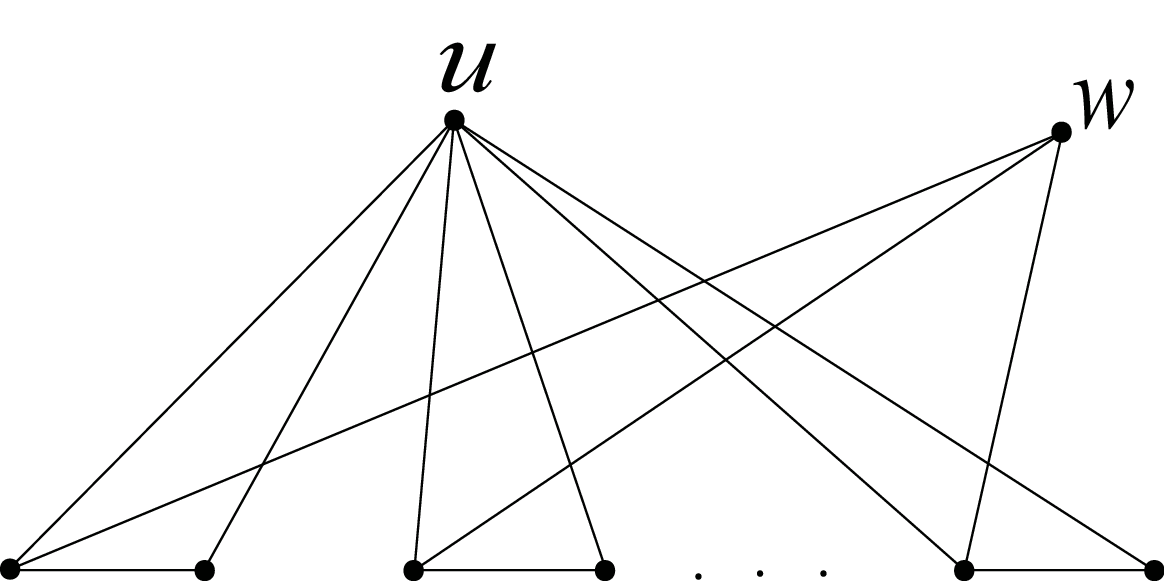}
         \caption{The graph $G_{2}$.}
         \label{2.}
   \end{minipage}
 \end{center}
\end{figure}

{\bf Subcase 2.2.} If $d(u)=n-2$, then by \eqref{2.}, we have
\begin{equation}\label{eq:ch-3}
\begin{split}
q(G^{\ast})&\leq d(u)+\frac{1}{d(u)}\sum_{w\in N(u)}d(w)\\&= n-t(\frac{n-3-d(u)}{d(u)})\\&= n+\frac{t}{n-2}.
\end{split}
\end{equation}

In this case, we consider two subcases as follows.

{\bf Subcase 2.2.1.} $t=0$.

By \eqref{eq:ch-3}, we have $q(G^{\ast})\leq n<q(F_{n})$, a contradiction.

{\bf Subcase 2.2.2.} $t=1$.

In this case, by Lemma \ref{le:ch-2.8}, we have $$n+\frac{1}{n-2}<n+\frac{2}{n-1}<q(F_{n})$$ for odd $n\geq6$ and $$n+\frac{1}{n-2}<n+\frac{2}{n}<q(F_{n})$$ for even $n\geq6$, a contradiction.

{\bf Subcase 2.2.3.} $2\leq t\leq\lfloor\frac{d(u)}{2}\rfloor$.

In this case, By \eqref{eq:ch-3}, either
$$
\begin{aligned}
q(G^{\ast})\leq n+\frac{t}{n-2}\leq n+\frac{1}{2}-\frac{1}{2(n-2)},
\end{aligned}
$$ with equality if and only if $e(N(u),N^{2}(u))=(n-1-d(u))(d(u)-t)$ and $d(u)=n-2, t=\frac{n-3}{2}$ for odd $n$ and hence $G^{\ast}\cong G_{1}$ (see Fig. \ref{1.}),
or
$$
\begin{aligned}
q(G^{\ast})\leq n+\frac{t}{n-2}\leq n+\frac{1}{2},
\end{aligned}
$$ with equality if and only if $e(N(u),N^{2}(u))=(n-1-d(u))(d(u)-t)$ and $d(u)=n-2, t=\frac{n-2}{2}$ for even $n$ and hence $G^{\ast}\cong G_{2}$ (see Fig. \ref{2.}). Obviously, neither $G_{1}$ or $G_{2}$ is not semi-regular bipartite graph or regular graph. Thus, $$2\leq t\leq\lfloor\frac{d(u)}{2}\rfloor-1$$ from Lemma \ref{le:ch-2.4}. We use similar above ways to consider $2\leq t\leq\lfloor\frac{d(u)}{2}\rfloor-1$ and there does not exist a extremal graph for $2\leq t\leq\lfloor\frac{d(u)}{2}\rfloor-1$. Actually, since $t\geq2$, we obtain that $G^{\ast}[N[u]]$ contains a $C_{3}$. Furthermore, it is impossible that $G^{\ast}$ is a bipartite graph. Since $d(u)=n-2$, we obtain that $G^{\ast}$ is not a regular graph. $\blacksquare$

\section{Proof of Theorem 1.3}

\noindent{\bf Proof.} Since $G^{\ast}$ is $\theta(1,2,3)$-free, we obtain that $G^{\ast}[N(u)]$ is $P_{4}$-free, and $e(N(u))\leq d(u)$ from Lemma \ref{le:ch-2.2}. Since $S_{n,2}$ is $\theta(1,2,3)$-free, combining with Lemma \ref{le:ch-2.9}, we have $n+2-\frac{4}{n+1}<q_{1}(S_{n,2})\leq q(G^{\ast})$. Let $u\in V(G^{\ast})$ be the vertex such that equality \eqref{eq:ch-1} holding.

If $\Delta(G^{\ast})=n-1$, then let $u\in V(G^{\ast})$ and $d(u)=n-1$. Thus, $u$ is a dominating vertex of $G^{\ast}$, and $4\leq e(N(u))\leq n-2$ from Lemma \ref{le:ch-2.6}. Furthermore, $q(G^{\ast})$ attains the maximum if and only if $G^{\ast}\cong S_{n,2}$. We will consider the two cases $e(N(u))\geq n-1$ and $e(N(u))\leq3$ as follows.

{\bf Case 1.} $e(N(u))\geq n-1$.

In this case, by Lemmas \ref{le:ch-2.4} and \ref{le:ch-2.5}, we have

$$n+2-\frac{4}{n+1}<q(G^{\ast})\leq d(u)+\frac{1}{d(u)}\sum_{w\in N(u)}d(w)\leq \frac{2e(G^{\ast})}{n-1}+n-2.$$
we get $$ e(N(u))> n-3+\frac{4}{n+1},$$ i.e., $$ e(N(u))\geq n-2.$$ Since $G^{\ast}[N(u)]$ is $P_{4}$-free, we have $e(N(u))\leq n-1$ from Lemma \ref{le:ch-2.3}. Combining above, $e(N(u))= n-1$ and $G^{\ast}[N(u)]\cong \frac{n-1}{3}K_{3}$ from Lemma \ref{le:ch-2.3}. Furthermore, we have $G^{\ast}\cong K_{1}\vee\frac{n-1}{3}K_{3}$. By carefully calculation, we obtain $$ q(K_{1}\vee\frac{n-1}{3}K_{3})=\frac{n+4+\sqrt{(n-4)^{2}+16}}{2}<q(S_{n,2}),$$ a contradiction.

{\bf Case 2.} $e(N(u))\leq 3$.

In this case, we obtain $G^{\ast}\subsetneq K_{1}\vee\frac{n-1}{3}K_{3}$ or $G^{\ast}\subsetneq S_{n,2}$, a contradiction.

We will consider $\Delta\leq n-2$ and we give three cases as follows.

{\bf Case 1.} If $d(u)=1$, by Lemma \ref{le:ch-2.4}, then we get
$$q(G^{\ast})\leq d(u)+\frac{1}{d(u)}\sum_{w\in N(u)}d(w)=1+\Delta\leq n-1<n+2-\frac{4}{n+1},$$ a
contradiction.

{\bf Case 2.} If $d(u)=2$, by Lemma \ref{le:ch-2.5}, then we get
$$q(G^{\ast})\leq d(u)+\frac{1}{d(u)}\sum_{w\in N(u)}d(w)=2+\Delta\leq n<n+2-\frac{4}{n+1},$$ a contradiction.

{\bf Case 3.} If $3\leq d(u)\leq n-2$, then we give the claim as follows.

{\bf Claim 4.2.} $e(N(u))\geq d(u)-1$.

{\bf Proof.} Let $u\in V(G^{\ast})$ be the vertex such that equality \eqref{eq:ch-1} holding. By Lemma \ref{le:ch-2.9}, we obtain
$$
\begin{aligned}
n+2-\frac{4}{n+1}<q(G^{\ast})&\leq d(u)+\frac{1}{d(u)}\sum_{w\in N(u)}d(w) \\&=d(u)+\frac{1}{d(u)}(d(u)+2e(N(u))+e(N(u),N^{2}(u)))\\&=d(u)+\frac{1}{d(u)}(d(u)+2e(N(u))+d(u)(n-1-d(u)),
\end{aligned}
$$
and hence $$e(N(u))>d(u)-\frac{2d(u)}{n+1}\geq d(u)-\frac{2(n-2)}{n+1}\geq d(u)-2,$$ due to $d(u)\leq n-2$ and $n\geq 6$. $\square$

By Claim 4.2, we obtain $G^{\ast}[N(u)]$ contains cycles or trees as subgraphs. Suppose that $G^{\ast}[N(u)]$ contains no cycles, we obtain that $G^{\ast}[N(u)]$ is a star from Claim 4.2. For any $w\in V(G^{\ast})\backslash N[u]$, we get $|N_{N(u)}(w)|\leq1$. Therefore, $e(N(u),N^{2}(u))\leq n-1-d(u)$. Moreover,
$$
\begin{aligned}
n+2-\frac{4}{n+1}<q(G^{\ast})&\leq d(u)+\frac{1}{d(u)}\sum_{w\in N(u)}d(w)\\&= d(u)+\frac{d(u)+2(d(u)-1)+n-1-d(u)}{d(u)}\\&=d(u)+2+\frac{n-3}{d(u)}\\&\leq max \{\frac{n+12}{3}, n+1-\frac{1}{n-2}\}\\&\leq n+1-\frac{1}{n-2}.
\end{aligned}
$$
Since $h(x)=x+\frac{n-3}{x}$ is a convex function for $x>0$, combining $n\geq6$ and $3\leq d(u)\leq n-2$, a contradiction. Hence, $G^{\ast}[N(u)]$ contains cycles as subgraphs and all cycles are $K_{3}$. Let $x,y,z $ be the vertices of a $K_{3}$. For any $w\in V(G^{\ast})\backslash N[u]$, we get $|N_{K_{3}}(w)|\leq1$. Therefore,
 $$ e(N(u),N^{2}(u))\leq (d(u)-2)(n-1-d(u)).$$
Moreover,
 $$
 \begin{aligned}
 n+2-\frac{4}{n+1}<q(G^{\ast})&\leq d(u)+\frac{1}{d(u)}\sum_{w\in N(u)}d(w)
 \\&= d(u)+\frac{d(u)+2e(N(u))+(d(u)-2)(n-1-d(u))}{d(u)}.
 \end{aligned}
 $$
 Furthermore,
 $$ e(N(u))>n-1-\frac{2d(u)}{n+1}>n-3+\frac{6}{n+1},$$ due to $d(u)\leq n-2$. In other words, $e(N(u))\geq n-2$. By Lemma \ref{le:ch-2.3}, we have $$ e(N(u))\leq d(u)\leq n-2. $$ Hence, $e(N(u))=d(u)=n-2$ and $G^{\ast}[N(u)]$ is a union of disjoint triangles. Furthermore, there only exists a vertex $w\in V(G^{\ast})\backslash N[u]$ and $|N_{N(u)}(w)|\leq1$. As above, we have $$ e(N(u),N^{2}(u))\leq 1.$$ Moreover,
 $$
\begin{aligned}
 n+2-\frac{4}{n+1}<q(G^{\ast})&\leq d(u)+\frac{1}{d(u)}\sum_{w\in N(u)}d(w)\\&= d(u)+\frac{d(u)+2d(u)+1}{d(u)}\\&=n+1+\frac{1}{n-2},
\end{aligned}
 $$
 which is a contradiction for $n\geq6$, as desired. $\blacksquare$

\section{Proof of Theorem 1.4}

\noindent{\bf Proof.} Since $G^{\ast}$ is $\{\theta(1,2,2),F_{5}\}$-free, we obtain that $G^{\ast}[N(u)]$ contains at most an edge and some isolated vertices.

\begin{lemma}\label{le:ch-6.1}
Let $q(S_{n,1}^{+})$ be the largest root of $x^{3}-(n+3)x^{2}+3nx-4=0$, and satisfy $$n<q(S_{n,1}^{+})<n+1 $$ for $n\geq 6$.
\end{lemma}
{\bf proof.} By Lemmas \ref{le:ch-2.1} and \ref{le:ch-2.2}, we get the quotient matrix of $Q(S_{n,1}^{+})$ is as follows:
$$
\left(\begin{array}{cccccc}
 n-1&2&n-3\\
 1&3&0\\
 1&0&1
\end{array}\right)
.$$
Then the characteristic polynomial of the quotient matrix of $Q(S_{n,1}^{+})$ is $$f(x)=x^{3}-(n+3)x^{2}+3nx-4.$$ Furthermore, $$f^{\prime}(x)=3x^{2}-2x(n+3)+3n>0$$ for $x\geq n$. By calculation, $f(n)<0$ and $f(n+1)>0$ for $n\geq6$, as desired. $\square$

Obviously, $S_{n,1}^{+}$ is $\{\theta(1,2,2),F_{5}\}$-free, we have $$n<q(S_{n,1}^{+})\leq q(G^{\ast}).$$ Let $u\in V(G^{\ast})$ be the vertex such that equality \eqref{eq:ch-1} holding.

If $\Delta(G^{\ast})=n-1$, then $G^{\ast}$ is a subgraph of $S_{n,1}^{+}$, by Lemma \ref{le:ch-2.7}, we have $q(G^{\ast})\leq q(S_{n,1}^{+})$, with equality if and only if $G^{\ast}\cong S_{n,1}^{+}$. We consider $\Delta(G^{\ast})\leq n-2$ as follows.

{\bf Case 1.} If $d(u)=1$, by Lemmas \ref{le:ch-2.4}, then we obtain $$n<q(G^{\ast})\leq d(u)+\frac{\sum_{uv\in E(G)} d(v)}{d(u)}=1+\Delta\leq n-1,$$ a contradiction.

{\bf Case 2.} If $2\leq d(u)\leq n-2$, then we consider two subcases as follows.

{\bf Subcase 2.1.} $G[N(u)]$ contains no an edge.

In this case, we have $e(N(u),N^{2}(u))\leq (n-1-d(u))d(u)$. By Lemmas \ref{le:ch-2.4}, we obtain
$$
\begin{aligned}
n<q(G^{\ast})&\leq d(u)+\frac{1}{d(u)}\sum_{w\in N(u)}d(w)\\&= d(u)+\frac{d(u)+2e(N(u))+e(N(u),N^{2}(u))}{d(u)}\\&= d(u)+\frac{d(u)+(n-1-d(u))d(u)}{d(u)}\\&=n,
\end{aligned}
$$
a contradiction.

{\bf Subcase 2.2.} $G^{\ast}[N(u)]$ contains an edge.

In this case, we have $e(N(u),N^{2}(u))\leq (n-1-d(u))(d(u)-1)$. By Lemma \ref{le:ch-2.4}, we obtain

\begin{equation}\label{eq:ch-4}
\begin{split}
q(G^{\ast})&\leq d(u)+\frac{1}{d(u)}\sum_{w\in N(u)}d(w)\\&=d(u)+\frac{d(u)+2e(N(u))+e(N(u),N^{2}(u))}{d(u)}\\&= d(u)+\frac{d(u)+2+(n-1-d(u))(d(u)-1)}{d(u)}\\&=n+1-\frac{n-3}{d(u)}.
\end{split}
\end{equation}

We consider the two subcases in the following.

{\bf Subcase 2.2.1.} If $2\leq d(u)\leq n-3$, by \eqref{eq:ch-4}, then we have $$q(G^{\ast})\leq n+1-\frac{n-3}{d(u)}\leq n<q(S_{n,1}^{+}),$$ a contradiction.

\begin{figure}[H]
 \begin{center}
    \includegraphics[height=4cm,width=5cm]{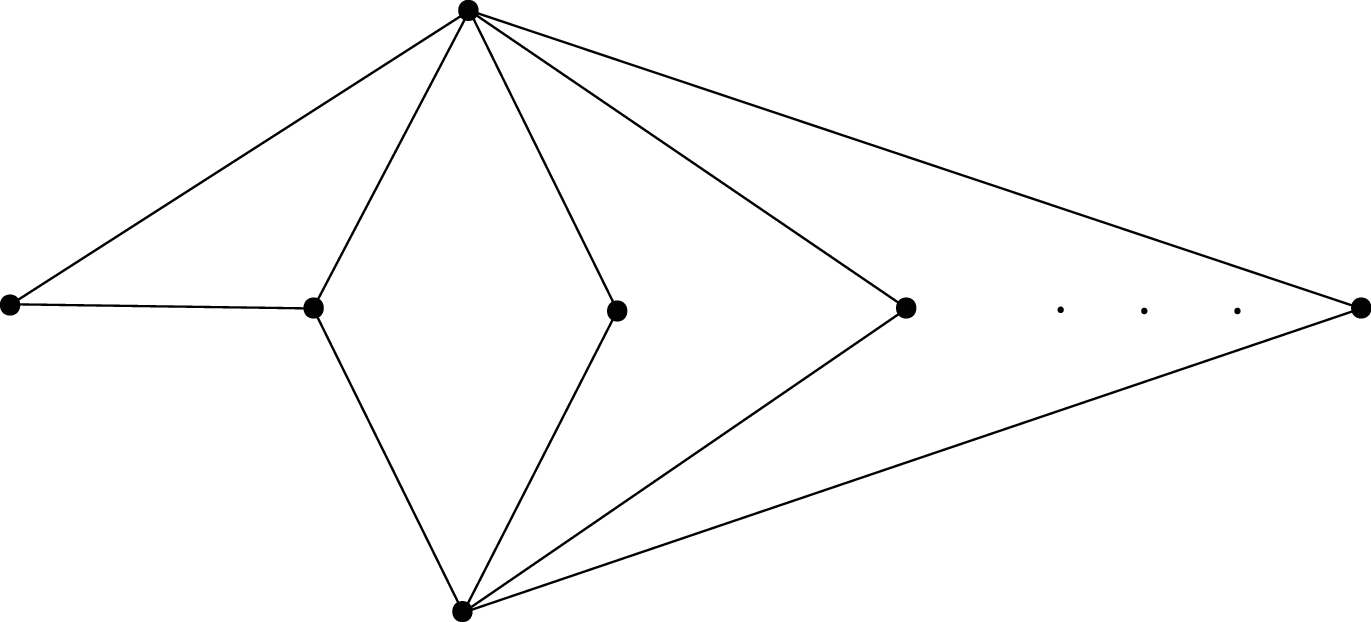}
     \caption{The graph $G_{3}$.}
      \label{3.}
      \end{center}
   \end{figure}

{\bf Subcase 2.2.2.} If $d(u)=n-2$, by \eqref{eq:ch-4}, we have $$q(G^{\ast})\leq n+1-\frac{n-3}{d(u)}\leq n+\frac{1}{n-2},$$ with equality $q(G^{\ast})=n+\frac{1}{n-2}$ if and only if $d(u)=n-2$ and $e(N(u),N^{2}(u))=(n-1-d(u))(d(u)-1)=n-3$. Thus, $G\cong G_{3}$ (see Fig. \ref{3.}). Since $G_{3}$ contains $C_{3}$ as a subgraph and $d(u)=n-2$, we obtain $G_{3}$ is neither a semi-regular bipartite graph or a regular graph. By Lemma \ref{le:ch-2.4}, $G\ncong G_{3}$ and $e(N(u),N^{2}(u))<(n-1-d(u))(d(u)-1)$, i.e., $ e(N(u),N^{2}(u))\leq (n-1-d(u))(d(u)-1)-1$. Moreover, we obtain
$$
\begin{aligned}
n<q(G^{\ast})&\leq d(u)+\frac{1}{d(u)}\sum_{w\in N(u)}d(w)\\&=d(u)+\frac{d(u)+2e(N(u))+e(N(u), N^{2}(u))}{d(u)}\\&= d(u)+\frac{d(u)+2+(n-1-d(u))(d(u)-1)-1}{d(u)}\\&=n+1-\frac{n-2}{d(u)}\\&\leq n,
\end{aligned}
$$
a contradiction.

This completes the proof. $\blacksquare$

\section*{Declaration of competing interest}

The authors declare that they have no conflict of interest.

\section*{Data availability}

No data was used for the research described in the article.

\end{document}